\documentclass[reqno]{amsart}
 
\usepackage{amsmath}
\usepackage{amsthm}
\usepackage{amssymb}
\usepackage{epsfig}
\usepackage{psfrag}
\usepackage{graphicx}
\usepackage{caption2}
\usepackage{textcomp}

\newcommand{\Rz}{\mathbb{R}}
\newcommand{\Nz}{\mathbb{N}}
\newcommand{\ffi}{\varphi}

\newcommand{\disp}{\displaystyle}

\newcommand{\ove}{\tilde}

 \newtheorem{theorem}{Theorem}[section]

\newtheorem{corollary}[theorem]{Corollary}

 \newtheorem{lemma}[theorem]{Lemma}
\begin{document}

\title[Counterexample to $\, C^{2,1}\,$ regularity]{A counterexample to $\, C^{2,1}\,$ regularity \\
for parabolic fully nonlinear equations} 

\author{Luis A. Caffarelli}
\address{Dept. of Mathematics, Univ. of Texas at Austin, Texas, USA.}
\email{caffarel\,@\,math.utexas.edu}
\urladdr{http://www.ma.utexas.edu/users/caffarel/}

\author{Ulisse Stefanelli}
\address{IMATI - CNR, 
v. Ferrata 1, I-27100 Pavia, Italy.}
\email{ulisse.stefanelli\,@\,imati.cnr.it} 
\urladdr{http://www.imati.cnr.it/ulisse/}
\thanks{Ulisse Stefanelli gratefully acknowledges the hospitality of the Department of Mathematics and the Institute for Computational Engineering and Sciences of the University of Texas at Austin, where most of this research was performed.}

\begin{abstract}
We address the self-similar solvability of a singular parabolic problem and show that solutions to parabolic fully nonlinear equations are not expected to be $\, C^{2,1}$.
\end{abstract}

\maketitle

\pagestyle{myheadings}

                                %
                                %
                                %
\section{Introduction}
\setcounter{equation}{0}

Locally bounded solutions to the fully nonlinear parabolic equation $\, u_t =F(D^2u)\,$ where $\, F \,$ is uniformly elliptic and convex are known to be locally $\, C^{2,\alpha}\,$ for some $\, \alpha \in (0,1)$ \cite{Bourgoing04,Wang0,Wang1,Wang2,Wang3}. The aim of this note is to show that $\, C^{2,1}\,$ regularity is generally not to be expected. To this end, we focus from the very beginning on the simplest fully nonlinear parabolic equation, namely
\begin{equation}\label{eq}
u_t = \max\big\{\Delta u, \Delta u/2\big\}.
\end{equation}
The latter equation has some interest in itself since it arises in a variety of different applicative frames ranging from control theory \cite{Bardi-Capuzzo97,Bensoussan-Lions84,LionsI,LionsII,LionsIII,Lions-Menaldi82}, to mechanics \cite{Barenblatt52,Barenblatt90}, combustion \cite{Buckmaster,Linan74}, biology \cite{Crooks04,Dancer99,Hilhorst03}, and finance \cite{Avellaneda,AP,L,HWW}. 

The $\, C^{2,\alpha}\,$ regularity theory for \eqref{eq} is due to {\sc Evans \& Lenhart} \cite{Evans-Lenhart81}. By letting  $\, \gamma (r):= (3 + \text{sign}(r))/4 \,$, equation \eqref{eq} can be rewritten as $\, u_t = \gamma(\Delta u) \Delta u\,$ or even $\, \beta(u_t)=\Delta u \,$ for $\, \beta (r) := 2r - r^+$. In particular, by differentiating \eqref{eq} with respect to time and letting $\, w = u_t\,$ we get that
\begin{equation}\label{eq2}
\beta(w)_t = \Delta w. 
\end{equation}

The H\"older regularity of $\, w \,$ follows from classic parabolic theory \cite{DiBenedetto93,LUS} by observing that $\,v=\beta(w)\,$ solves the equation in divergence form $\, v_t = \text{\,div\,} ( \gamma(v)\nabla v)$. 

We shall prove here that $\, w \,$ cannot be expected to be Lipschitz continuous, which implies that the corresponding $\, u \,$ fails to be $\, C^{2,1}$. To this aim, we will provide a full description of radial self-similar solutions to \eqref{eq2}. In particular, we are here interested in locally bounded
positively-homogeneous solutions, namely 
$$\, w (x,t) = (-t)^{\alpha/2}f(|x|/{ \sqrt{-t}})\quad  x \in \Rz^N , \ t <0, \ \alpha>0$$
for some smooth self-similar profile $\,f $. We prove that, in order for a locally bounded self-similar solution to exist, some specific choice for $\, \alpha \,$ has to be made.

Our main result reads as follows.

\begin{theorem}\label{nuovomain} There exist two sequences $\,\alpha_k^+\,$ and $\, \alpha_k^- \,$ such that 
\begin{enumerate}
\item[i)] a locally bounded $\, \alpha$-homogeneous solution $\,w \,$ with $\, w(0,-1)=\pm 1\,$ exists iff $\, \alpha = \alpha_k^\pm \,$ for some $\,k \in \Nz$,
\item[ii)] $\,\alpha_k^+\,$ and $\, \alpha_k^- \,$ are such that the corresponding self-similar profile changes sign exactly $\, k \,$ times,
\item[iii)] $\, \alpha^\pm_k \,$ are strictly increasing and unbounded,
\item[iv)] $\,0< \alpha^-_1 < 2 < \alpha^+_1$,
\item[v)] $\,\alpha^-_k< \alpha^+_{k+1}<\alpha^-_{k+2} \ $ for all $\,  k \in \Nz$.
  \end{enumerate}
\end{theorem}

In particular, we prove that there exists a locally bounded solution $\, w \,$ to \eqref{eq2} such that $\, w(0,t)=-(-t)^{\alpha^-_1/2} \,$ for $\, t \in (-1,0)$. Hence, owing to Theorem \ref{nuovomain}.iv, one has that $\, w \,$ fails to be Lipschitz and we have the following.

\begin{corollary}\label{c21no}
  Equation \eqref{eq} admits locally bounded solutions which are not $\, C^{2,1}$.
\end{corollary}

A remarkable feature of this construction is that it is not symmetric with respect to sign changes. Namely, the only (up to multiplication by positive constants) non-Lipschitz locally bounded radial self-similar solution $\, w \,$ to \eqref{eq2} is such that $\, w(0,-1) < 0 $ and $\, w(\cdot,-1) \,$ has exactly one sign change, i.e. $\, w(\cdot,-1) \,$ is negative inside some given ball centered in $\, 0 \,$ and positive outside. In particular, all locally bounded radial self-similar solutions $\, w \,$ to \eqref{eq2} with $\, w(0,-1) >0 \,$ turn out to be Lipschitz. The same Lipschitz regularity holds for all locally bounded radial self-similar solutions such that $\, w(\cdot,-1) \,$ changes sign more than once.

Clearly, the counterexample to $\, C^{2,1}\,$ regularity for \eqref{eq} is not requiring the full extent of Theorem \ref{nuovomain}. We however provide here a detailed description of the self-similar solvability of \eqref{eq2} for the sake of completeness. In particular, we directly focus on solutions for $\, \alpha >0 \,$ in Sections \ref{selfsimilar_solutions}-\ref{proof} and then, by developing some suitable duality technique based of the so called Appell transformation, we discuss the case of self-similar solutions for $\, \alpha<0 \,$ in Section \ref{Appell}.

Let us recall that the related case of solutions of \eqref{eq} of negative homogeneity has been already considered by {\sc Kamin, Peletier, \& Vazquez} \cite{Aronson-Vazquez95,Kamin-Peletier-Vazquez89,Kamin-Peletier-Vazquez91}. Our analysis, although completely independent, presents some reminiscence of the former as we comment in Section \ref{Appell} below.

\section{Self-similar solutions}\label{selfsimilar_solutions}
\setcounter{equation}{0}

We are interested in positively-homogeneous solutions of \eqref{eq2}, namely functions $\, w \,$ such that $\, w(\mu x, \mu^2 t) = \mu^{\alpha}w(x,t) \,$ for all $\, \mu>0, \, x \in \Rz^N, \, t <0\,$ and some $\, \alpha > 0$. We are thus imposing
$$\, w (x,t) = (-t)^{\alpha/2}\ffi(x/{ \sqrt{-t}})\quad x\in \Rz^N, \ t <0$$
for some smooth function $\, \ffi$. 

Plugging the latter expression into the relation $\, \lambda w_t = \Delta w \,$ for $\, \lambda >0\,$ one easily checks that $\, \ffi \,$ solves
$$- \Delta \ffi +\frac{\lambda}{2}\frac{x}{\sqrt{-t}} \cdot \nabla \ffi = \frac{\lambda \alpha}{2} \ffi.$$
Let us assume now $\, \ffi \,$ to be radially symmetric and use the symbol $\, s \,$ for the self-similar variable $\,  |x|/ \sqrt{-t}$. Hence, the function $\, f(s) = \ffi(x/ \sqrt{-t})\,$ fulfills 
\begin{equation}\label{eq_ffi}
 f''(s) + \left( \frac{N-1}{s} -\frac{\lambda s}{2} \right)f'(s) + \frac{\lambda \alpha }{2}f(s)=0 .
\end{equation}
In view of the expected regularity of $\, w \,$ at the origin for $\, t <0$, we must have $\, f'(0)=0$. On the other hand, since we aim at locally bounded solutions $\, w $, we will require $\, s \mapsto s^{-\alpha}f(s) \,$ to remain bounded as $\, s\,$ goes to $\, +\infty$.
 
In particular, $\, f \in C^{1,1}[0,+\infty)\,$ is said to be a self-similar profile if
\begin{eqnarray}
\quad\left\{
\begin{array}{l}
f(0)\not = 0, \quad f'(0)=0,\quad\\
\disp  f''(s) + \left( \frac{N-1}{s} -s \right)f'(s) + \alpha f(s)=0 \quad \text{where} \ \ f <0, \ \ s>0 \\
\disp f''(s) + \left( \frac{N-1}{s} -\frac{ s}{2} \right)f'(s) + \frac{ \alpha }{2} f(s)=0\quad \text{where} \ \ f >0, \ \ s>0 .
\end{array}
\right.\label{S}
\end{eqnarray}
and, additionally,
\begin{equation}
  \quad  s \mapsto s^{-\alpha}f(s) \ \  \text{is bounded.}\label{bound}
\end{equation}

Clearly, Problem \eqref{S} can be solved for any $\, \alpha >0 \,$ by means of a direct construction argument. One has to find a solution to the corresponding Cauchy problem up to the first zero (which always exists, see Lemma \ref{zero} below). Then, one restarts a second Cauchy problem with the other equation. At this point, either the solution is unbounded or it has a zero. In this second case, one again solves a Cauchy problem and the procedure goes on up to finite termination (see below). 

Among all candidate self-similar profiles constructed as above, we prove that the boundedness assumption \eqref{bound} can be fulfilled only by specific choices of $\,\alpha$. In particular, Theorem \ref{nuovomain} entails that Problem \eqref{S}-\eqref{bound} admits countably many solutions, each of which is uniquely determined (up to multiplication by positive constants) by the sign of $\, f(0) \,$ and the number of sign changes. All other self-similar profiles will show exponential growth at infinity instead.

By directly considering \eqref{eq2} or computing on \eqref{S}, we shall explicitly observe that 
\begin{gather}
 f>0  \ \ \text{solves} \ \  \eqref{S} \ \ \text{in the interval}\ \  I \subset [0,+\infty)  \nonumber\\
  \text{iff} \ \  g(s)= - f(\sqrt2 s)  \ \ \text{solves} \ \  \eqref{S} \ \ \text{in} \ \ I/\sqrt2.
\label{cambio}
\end{gather}
This change of variables will turn out to be useful later on.

Let us now draw a relation between Problem \eqref{S}-\eqref{bound} and some suitable eigenvalue problem. Namely, we remark that \eqref{eq_ffi} may be rewritten as
\begin{equation}\label{calcolino}
  - \big( s^{N-1}e^{-\lambda s^2 /4}f'(s) \big)' =s^{N-1} e^{-\lambda s^2 /4}\frac{\lambda \alpha}{2}f(s).
\end{equation}
Let now $\,\mu^+ \,$ and $\, \mu^-\,$ be the absolutely continuous measures given by 
$$\, d \mu^+(s) = s^{N-1}e^{-s^2/4}ds, \quad d\mu^-(s)=  s^{N-1}e^{-s^2/2} ds.$$ 
Moreover, for all intervals $\, I \subset [0,+\infty)$, we define the spaces
\begin{gather}
 H^{+}(I) = L^2(I,\mu^{+}), \quad V^+(I)=\{f \in H^{+}(I) \ : \ f' \in H^{+}(I) \}\nonumber\\
W^+(I)= \{f \in V^+(I)\ : \ f(\iota) =0 \ \ \text{if} \ \ \iota \in \{\inf I, \sup I\}\setminus \{0, +\infty\}\},\nonumber
\end{gather}
and the spaces $\, H^-(I), \, V^-(I)$, and $\, W^-(I) \,$ correspondingly by means of the measure $\, \mu^-$. Let us consider the eigenvalue problems
\begin{eqnarray}
&&\text{find} \quad f \in W^+(I), \ f\not = 0, \ \ \text{such that}\label{A^+}\\
&& \qquad\qquad \int_I f'g' \, d\mu^+=\disp\frac{\alpha}{2}\int_I f g \, d\mu^+\quad \forall g \in W^+(I)\nonumber\\
&&\text{find} \quad f \in W^-(I), \ f \not =0,   \ \ \text{such that}\label{A^-}\\
&&  \qquad\qquad \int_I f'g' \, d\mu^-=\disp{\alpha}\int_I f g \, d\mu^-\quad \forall g \in W^-(I)\nonumber
\end{eqnarray}
It is straightforward to check that a non-negative $\, f \in W^+(I) \,$ solves \eqref{calcolino} iff it solves \eqref{A^+}. On the other hand, a non-positive $\, f \in W^-(I) \,$ solves \eqref{calcolino} iff it solves \eqref{A^-}. 
For the sake of later purposes, let us introduce
some notation for the Rayleigh quotients
$$R^+(f,I)= \disp\frac{\disp\int_I (f')^2 d\mu^+}{\disp\int_I f^2 d\mu^+}, \quad R^-(f,I)= \disp\frac{\disp\int_I (f')^2 d\mu^-}{\disp\int_I f^2 d\mu^-},$$
which are defined for $\, f \not = 0 \,$ and $\, f \in V^+(I)\,$ ($ f \in V^-(I)$, respectively).  

\section{Proof of Theorem \ref{nuovomain}}\label{proof}
\setcounter{equation}{0}

Let us start by observing that the strictly positive lower bound in Theorem \ref{nuovomain}.iii is obviously ensued from the above recalled $\, C^{2,\alpha}\,$ regularity theory. Namely, all locally bounded solutions to \eqref{eq2} are $\, C^{\alpha/2}\,$ in time for some $\, \alpha >0\,$ depending just on $\, N \,$ (and on $\, \beta$). The proof of Theorem \ref{nuovomain} will follow from a direct construction argument.

\subsection{Shooting from $0$} Let us start from the following.

\begin{lemma}\label{zero}
For all $\, \alpha >0$, the solution $\, f \,$ to \eqref{S} has a zero.
\end{lemma}
\begin{proof}
Let us consider $\, f(0)>0\,$ (the case $\, f(0)<0\,$ being completely analogous).
 Taking the limit as $\, s \rightarrow 0^+ \,$ in the equation one gets that 
$$  \lim_{s\rightarrow 0^+}f''(s) = -\alpha f(0)/(2N) <0.$$
Hence $\, f' <0 \,$ at least locally in a right neighborhood of zero. Indeed, $\, f' <0 \,$ as long as $\, f >0$. By contradiction let $\, \bar s \,$ be the first point where $\, f' \,$ vanishes and assume $\, f>0 \,$ on $\, (0,\bar s)$. One computes from $\, (s^{N-1}f'(s))'=s^{N-1}(s f'(s) - \alpha f(s))/2\,$ that
$$0=\int_0^{\bar s}(s^{N-1}f'(s))'ds = \frac12\int_0^{\bar s} s^Nf'(s)\,ds - \frac{\alpha}{2}\int_0^{\bar s} s^{N-1}f(s)\,ds < 0, $$
a contradiction. Assume now $\, f >0 \,$ everywhere. Then, for $\, s^2 >2(N-1)\,$ one has that $\, f''(s) < -\alpha f(s)/2 <0\,$ and $\, f \,$ is concave and decreasing. Hence, it has a zero.
\end{proof}

Owing to the latter lemma, we will denote by $\, s^{+}_{\alpha}\,$ the first zero of $\, f^+_{\alpha} \,$ where $\,\alpha>0 \,$ is given and $\, f^+_{\alpha}\,$ solves \eqref{S} with $\, f^+_{\alpha}(0)>0\,$ (recall that the latter functions are defined in all of $\, [0,+\infty)$). One has that the unique solvability of \eqref{S} entails in particular that $\, (f^{+}_{\alpha})'(s^{+}_{\alpha})< 0$. At the same time, owing to \eqref{cambio}, we shall let $\, s^-_\alpha = s^+_\alpha/\sqrt2\,$ and notice that $\, s^-_\alpha \,$ is the first zero of a solution to \eqref{S} starting from a negative value.

We can now make precise the argument of Section \ref{selfsimilar_solutions} by observing the following. 
\begin{lemma}\label{polp}
  For any $\, \alpha >0$, the minimal eigenvalue of Problem \eqref{A^+} on $\, (0, s^{+}_{\alpha}) \,$ is $\, \alpha\,$ and the corresponding eigenfunction is $\, f^{+}_{\alpha}$.
\end{lemma}

\begin{proof}
Classical results \cite[Thm 8.38, p. 214]{Gilbarg-Trudinger} ensure that, for all $\,  s>0$, the minimal eigenvalue of Problem \eqref{A^+} on $\, (0, s) \,$ is strictly positive and simple and the corresponding eigenfunction is the only eigenfunction which does not change sign on $\, (0,s)$. On the other hand, $\, f^+_\alpha \,$ solves \eqref{S} and does not change sign on $\, (0,s^+_\alpha)$, whence the assertion follows.
\end{proof}

\begin{figure}
\centering
       \includegraphics[width= .45  \textwidth]{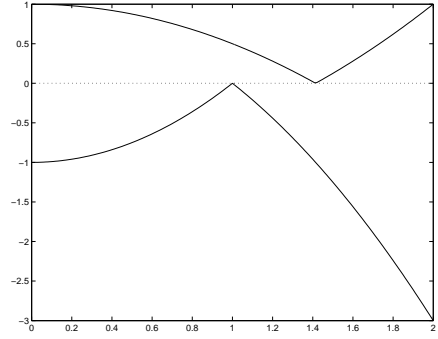}%
  \caption{The parabolas $\, |p^+|\,$ and $\, -|p^-|\ \ (N=1)$}
  \label{parabola}
\end{figure}

Clearly, the dual statements on Problem \eqref{A^-} on $\, (0, s^{-}_{\alpha}) \,$ hold true as well. By exploiting continuous dependence and the characterization of Lemma \ref{polp}, we easily deduce that the functions $\, \alpha \mapsto s^+_\alpha\,$ and $\, \alpha \mapsto s^-_\alpha \,$ are continuous, strictly decreasing, and onto $\, (0,+\infty)$. We shall denote their inverses by $\, \alpha^{+}\,$ and $\, \alpha^-$, respectively. Moving from \eqref{cambio}, we readily have that $\, \alpha^+(s) = \alpha^-(s/\sqrt2)\,$ for all $\, s \geq 0$. It is moreover a standard matter to consider the parabolas (see Figure \ref{parabola})
$$ p^+(s)= -s^2/(2N)+1\ \ \ \text{and} \ \ \ p^-(s)=- p^+(\sqrt2 s)= s^2/N -1\,$$
in order to check that $\,  s^+_2=\sqrt{2N}\,$ and $\,  s^-_2=\sqrt{N}$, or, equivalently
\begin{equation}\label{par1}
\alpha^+(\sqrt{2N})=2=\alpha^-(\sqrt{N}).
\end{equation}

Hence, let us summarize the above results by stating the following.

\begin{lemma}\label{sette}
The functions $\, \alpha^{+}\,$ and $\, \alpha^-\,$  are continuous, strictly decreasing, and onto.
\end{lemma}

Before closing this subsection, we shall explicitly remark that all solutions $\, f \,$ to \eqref{S} turn out to be definitely strictly monotone. Indeed, this may be proved by simply adapting the proof of Lemma \ref{zero} since a self-similar profile is forced to change sign after any critical point and all self-similar profiles change sign a finite number of times.

\subsection{Shooting from $+\infty$} Let us now consider the possibility of solving \eqref{S} for $\, f \,$ after its first zero (which always exists due to Lemma \ref{zero}). As mentioned above, two situations may occur. It could happen that $\, f\,$ never changes sign or that it comes back to zero again. In any case, by gluing together pieces, one can construct a $\, C^{1,1}\,$ solution to the differential problem in \eqref{S}.

We shall hence discuss the possibility of fulfilling the bound \eqref{bound} by proving the following.
\begin{lemma}\label{plus}
For any $\, \alpha >1\,$ there exists $\, \ove s^+_\alpha >0 \,$ such that,
letting $\, f\,$ solve \eqref{S} for  $\, s>  \ove s^+_\alpha
\,$ with $\,  f  (\ove s^+_\alpha)=0 <  f'  (\ove s^+_\alpha)$, one has that $\, f \,$ is positive for  $\,s> \ove s^+_\alpha\,$ and  $\, s \mapsto s^{-\alpha} f(s)\,$ is bounded.
\end{lemma}
\begin{proof} We change variables by letting $\, h(t):= t^\alpha f(1/t)$, and transform \eqref{S} into
\begin{eqnarray}
\left\{
\begin{array}{l}
h(0)\not = 0, \quad h'(0)=0,
\\
\disp  t^2h''(t) - \left({2\alpha} + N - 3  \right) t h'(t) + \frac{h'(t)}{t}   + \alpha(\alpha + N - 2) h(t) =0 \\
 \qquad \text{if} \ \ h <0, \\
\disp  t^2h''(t) - \left({2\alpha} + N - 3  \right) t h'(t) + \frac{h'(t)}{2t}  + \alpha(\alpha + N - 2) h(t) =0\\
\qquad    \text{if} \ \ h >0.
\end{array} 
\right.\label{S2}
\end{eqnarray}
The initial condition $\, h'(0)=0 \,$ is forced by the fact that we ask the latter singular differential problem to be uniquely integrable. 

We shall check that, for any $\, \alpha >1$, the solution $\, h \,$ to \eqref{S2} admits indeed a zero. Let us focus from the very beginning on the case $\, h(0)>0\,$ (the other case being analogous) and check that 
$$ \lim_{t \rightarrow 0^+} h''(t) = - 2\alpha(\alpha +N - 2)h(0)<0.$$
 Namely, $\, h' <0 \,$ in a right-neighborhood of $\, 0$. We have that $\, h' <0 \,$ as long as $\, h >0$. Indeed, $\, h \,$ cannot have a positive minimum in $\, t >0 \,$ since, in that case,
$$t^2h''(t) + \alpha(\alpha +N - 2) h(t) =0,$$
which leads to $\, h''(t) <0$, a contradiction. On the other hand, we have that
$\, -(2\alpha + N - 3) + 1/(2t) <0 \,$ for $\, t \,$ large and hence
$$t^2h''(t) + \alpha(\alpha +N - 2) h(t) <0\quad \text{for $\, t \,$ large},$$
and $\, h \,$ is concave and decreasing as long as $\, h>0$. Hence, it has a zero.

We denote by $\, t^+_\alpha >0 \,$ the first zero constructed above, changing variables back as $\, f^+_\alpha(s): = s^\alpha h(1/s)$, and letting $\, \ove s^+_\alpha = 1/ t^+_\alpha$, we have found the unique solution to \eqref{S} for  $\, s >  \ove s^+_\alpha\,$ with $\,  f^+_\alpha  (\ove s^+_\alpha)=0 <  ( f^+_\alpha)'  (\ove s^+_\alpha)$. The latter does not change sign and grows as $\, s^\alpha\,$ at infinity
\end{proof}

The restriction $\, \alpha >1 \,$ in Lemma \ref{plus} is actually needed in the latter proof only for $\, N = 1\,$ and one could ask $\, \alpha >1/2 \,$ for $\, N=2 \,$ and $\, \alpha >0 \,$ for $\, N \geq 3$, instead. On the other hand, as it will be clear later, the above mentioned restriction is irrelevant for the purposes of our analysis.

The self-similar profile $\, f^+_\alpha \,$ constructed in the proof of Lemma \ref{plus} is hence an eigenfunction of Problem \ref{A^+} on $\, (\ove s_\alpha^+,+\infty) \,$ and does not change sign. Again, this amounts to say that $\, \alpha \,$ is the minimal eigenvalue of Problem \ref{A^+} on $\, (\ove s_\alpha^+,+\infty) $. Namely, we have the following.

\begin{lemma}\label{polp2}
  For any $\, \alpha >1$, the minimal eigenvalue of Problem \eqref{A^+} on $\, (\ove s^{+}_{\alpha}, +\infty) \,$ is $\, \alpha\,$ and the corresponding eigenfunction is $\, f^{+}_{\alpha}$.
\end{lemma}

Clearly, again exploiting \eqref{cambio}, Lemmas \ref{plus} and \ref{polp2} can be rephrased for negative self-similar profiles as well. In particular, for all $\, \alpha >1\,$ one finds $\, \ove s^-_\alpha >0 \,$ such that any self-similar profile $\, f^-_\alpha\,$ fulfilling \eqref{S} for  $\, s>  \ove s^-_\alpha
\,$ with $\,  f^-_\alpha  (\ove s^-_\alpha)=0 >  ( f^-_\alpha)'  (\ove s^-_\alpha)\,$ does not change sign and is such that $\, s \mapsto s^{-\alpha} f^-_\alpha(s)\,$ is bounded. Moreover, $\, \alpha \,$ is the minimal eigenvalue of Problem \eqref{A^-} on $\, (\ove s^-_\alpha,+\infty)$. In particular, we have defined the functions $\, \alpha \mapsto \ove s^+_\alpha \,$ and  $\, \alpha \mapsto \ove s^-_\alpha\,$ which turn out to be strictly increasing. We shall denote their corresponding inverses by $\, \ove \alpha^+\,$ and $\, \ove \alpha^-$, respectively (recall that $\, \ove \alpha^+(s) = \ove \alpha^-(s/\sqrt2)\,$ for all $\, s>0$). By considering again the parabolas $\, p^+\,$ and $\, p^-\,$ we observe that
\begin{equation}\label{par2}
\ove \alpha^+(\sqrt{2N})=2= \ove \alpha^-(\sqrt{N}).
\end{equation}

Let us recall that we have checked the following.

\begin{lemma}\label{otto} The functions $\, \ove \alpha^+\,$ and $\,\ove \alpha^-\,$ are continuous and strictly increasing.\end{lemma}

In order to clarify the meaning of the functions $\, \ove \alpha^+\,$ and $\,\ove \alpha^-$ one can again exploit the relation with the eigenvalue Problems \eqref{A^+} and \eqref{A^-}. In particular, let us fix $\,\alpha >1 \,$ and the related $\,\ove s^+_\alpha$. We shall consider Problem \eqref{A^+} in a sequence of increasing intervals $\,I^r:= (\ove s^+_\alpha,r)$ for $\, r \to +\infty$. The corresponding minimal eigenvalues $\, \alpha^+_r\,$ are positive and simple and their eigenfunctions $\, f^+_r\,$ do not change sign on $\, I_r$. Moreover, the function $\,r \mapsto \alpha^+_r \,$ is continuous, strictly decreasing, and unbounded as $\, r \searrow s^+_\alpha$. By exploiting Lemma \ref{polp2} and standard compactness arguments, one can easily check that 
$$ \lim_{r\to +\infty} \alpha^+_r =\alpha,$$
and the corresponding eigenfunctions converge locally uniformly to the self-similar profile $\, f^+_\alpha$. Indeed, we have that $\, \lim_{r\to +\infty} \alpha^+_r =: \overline \alpha \,$ exists and, by suitably renormalizing and extracting not relabeled subsequences, there exists $\, f \,$ such that 
$$f_r \to f \quad \text{weakly in} \ \ W^+(\ove s^+_\alpha, +\infty).$$
Hence, the positivity of $\, f \,$ follows along with the fact that $\, f \,$ solves \eqref{S} with $\, \alpha \,$ replaced by $\, \overline \alpha$. Finally, $\, f \,$ is the eigenfunction corresponding to $\, \alpha \,$ by Lemma \eqref{polp2}, namely $\, \overline \alpha = \alpha$. Let us rephrase these facts in the following lemma.

\begin{lemma}\label{cosa}
  Let $\, \alpha,\, \beta > 1 \,$ be given and $\, f \,$ solve \eqref{S} with $\, \alpha \,$ replaced by $\, \beta \,$ for $\, s > \ove s^+_\alpha \,$ with $\,f(\ove s^+_\alpha) =0 <   f'(\ove s^+_\alpha)$. Then 
  \begin{enumerate}
 \item[i)] if $\, \alpha < \beta \,$ then $\, f \,$ changes sign,
\item[ii)] if $\, \alpha = \beta \,$ then $\, f \,$ does not change sign and grows as $\, s^\alpha\,$ at infinity,
 \item[iii)] if $\, \alpha > \beta \,$ then $\, f \,$ does not change sign and grows exponentially at infinity.
  \end{enumerate}
\end{lemma}
\begin{proof}
  Lemma \ref{plus} entails ii). As for i) and iii), owing to the above discussion, we readily check $\, f \,$ has a zero in $\, r \,$ iff $\, \beta = \alpha^+_r >\alpha\,$ and we have proved i). If $\, \beta < \alpha \,$ then surely $\, f \not \in W^+(\ove s^+_\alpha,+\infty)\,$ owing to the minimality of $\,\alpha$. In particular, the exponential growth of $\, f \,$ follows.
\end{proof}

\subsection{Zeros}
Let us collect here some remark from the above shooting constructions. Given any $\,\alpha>1$, let $\, f^\pm_\alpha \,$ denote the solutions to \eqref{S} starting from $\,f^\pm_\alpha(0)=\pm1\,$ and solving successively the corresponding Cauchy problems up to the last zero. We shall denote by $\, s^{\pm,k}_\alpha \,$ the $\, k$-th zero of $\, f^\pm_\alpha\,$ (namely, $\, k \,$ ranges on a finite set of indices which will be proved to be increasing with $\, \alpha$). Moreover, let $\, g^\pm_\alpha \,$ be the outcome of the variable transformation $\, s=1/t \,$ applied to the solution $\, h^\pm_\alpha \,$ of \eqref{S2} starting from $\,h^\pm_\alpha(0)=\pm1\,$ and built by solving successively the corresponding Cauchy problems up to the last zero. We let $\, \ove s^{\pm,k}_\alpha \,$ the $\, k$-th zero of $\, g^\pm_\alpha$, ordered starting from $\, s=0$.

The characterization  of $\, \alpha \,$ as minimal eigenvalue of Lemmas \eqref{polp} and \eqref{polp2} and the corresponding monotonicity properties of Lemmas \ref{sette} and \ref{otto} entail the following monotonicity properties for the zeros of $\, f^\pm_\alpha\,$ and $\, g^\pm_\alpha $.

\begin{lemma}\label{zz} The functions $\, \alpha\mapsto s^{\pm,k}_\alpha\,$ are continuous and strictly decreasing. The functions $\, \alpha\mapsto \ove s^{\pm,k}_\alpha\,$ are continuous and strictly increasing.  
\end{lemma}

\begin{proof} Continuity obviously follows from continuous dependence and the monotonicity of $\,\alpha\mapsto s^{\pm,1}_\alpha\,$ and $\, \alpha\mapsto \ove s^{\pm,k}_\alpha\,$ has already been proved in Lemmas \ref{sette} and \ref{otto}. Assume by contradiction that $\, s^{+,2}_\alpha\,$ is not strictly decreasing with $\, \alpha$. Namely, assume there exist $\,\alpha_1 < \alpha_2\,$ such that the corresponding self-similar profiles have at least two zeros and $\,s^{+,2}_{\alpha_1}\leq s^{+,2}_{\alpha_2}$. Letting $\, I_j:=(s^{+,1}_{\alpha_j},s^{+,2}_{\alpha_j})\,$ for $\, j=1,2\,$ and observing that $\, I_1 \subset I_2$, one has that
$$\alpha_2 = R^-(f^+_{\alpha_2},I_2) \leq R^-(f^+_{\alpha_1},I_1)=\alpha_1 < \alpha_2, $$
a contradiction. Once we have established that indeed $\,\alpha\mapsto s^{+,2}_\alpha\,$ is strictly decreasing, this same argument can be iterated inductively up to the last zero of $\, f^+_\alpha$. Clearly the monotonicity of $\, \alpha \mapsto s^{-,k}_\alpha \,$ follows with minor changes. Moreover, the very same proof can be exploited for $\, g^\pm_\alpha\,$ as well (from its last to its first zero though).
\end{proof}

\subsection{Construction of solutions changing sign once}

We shall now combine the results of the previous subsections in order to build self-similar profiles changing sign exactly once. By considering again the parabolas $\, p^+(s)= -s^2/(2N)+1\,$ and $\, p^-(s)= s^2/N -1\,$ and, in particular, owing to \eqref{par1} and \eqref{par2}, we readily check that we cannot find a solution to Problem \eqref{S} by gluing together two semi-profiles of homogeneity $\, 2\,$  (see Figure \ref{parabola}). In particular, we have observed that
$$ s^{+,1}_2 = \sqrt{2N} > \ove s^{-,1}_2=\sqrt{N}.$$
By means of Lemma \ref{zz}, we shall increase $\,\alpha\,$ and find the unique $\, \alpha^+_1 \,$ such that 
$$ s^{+,1}_{\alpha^+_1} =  \ove s^{-,1}_{\alpha^+_1}.$$
The two semi-profiles corresponding to this specific homogeneity $\, \alpha^+_1\,$ can be glued together in order to form a solution of \eqref{S}. Moreover, we readily have that $\, \alpha^+_1 >2\,$ by construction. 

On the other hand, starting from the fact that 
$$ s^{-,1}_2 = \sqrt{N} < \ove s^{+,1}_2=\sqrt{2N},$$
we argue as above and decrease $\, \alpha\,$ in order to check that there exists a unique $\, \alpha^-_1<2 \,$ such that 
$$ s^{-,1}_{\alpha^-_1} = \ove s^{+,1}_{\alpha^-_1}.$$
\begin{figure}
   \begin{minipage}[t]{0.5\linewidth}
      \centering
      \includegraphics[width= .9 \textwidth]{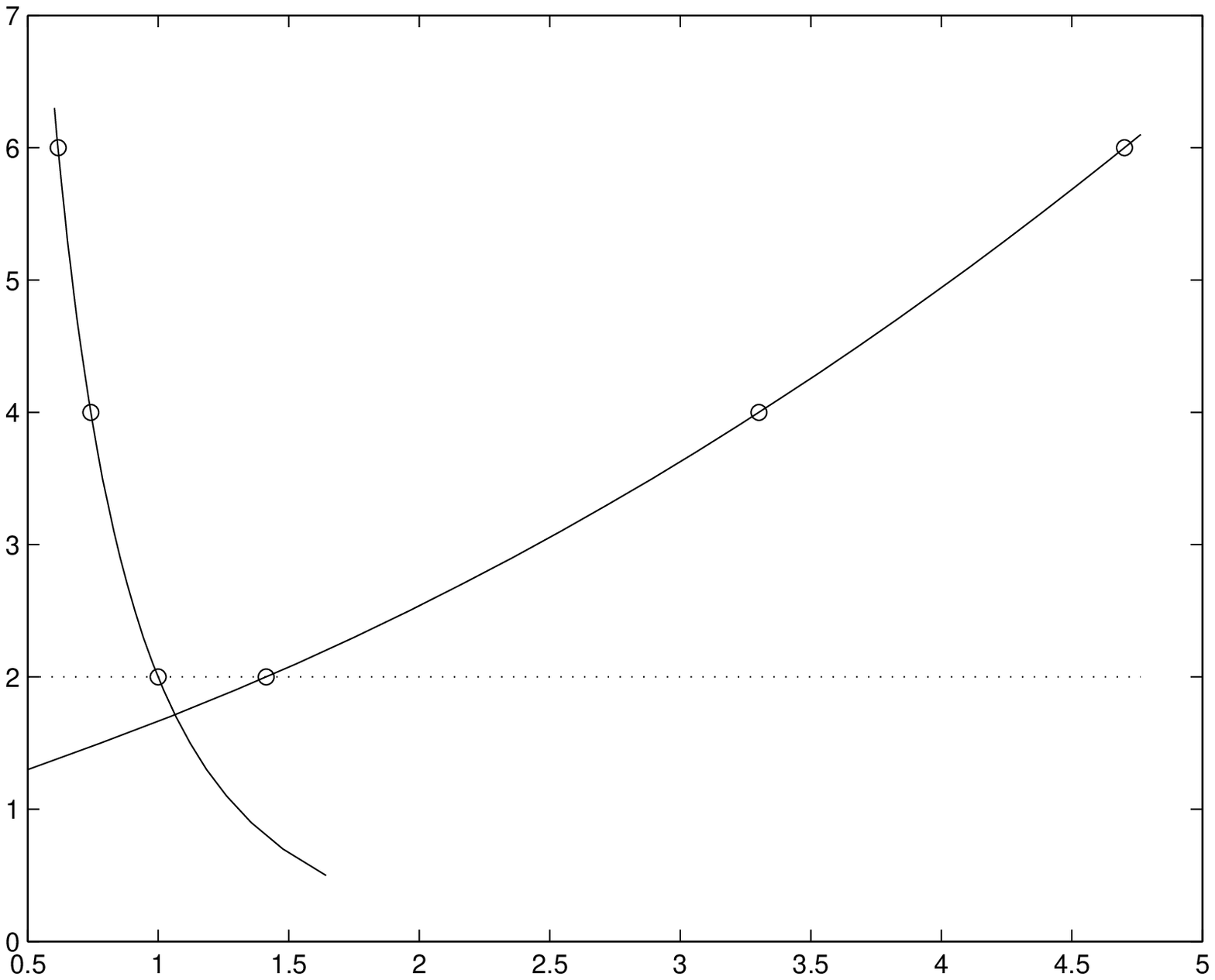}%
  \caption{Constructing $\, \alpha^-_1$.}
  \label{contro}
   \end{minipage}%
   \begin{minipage}[t]{0.5\linewidth} 
     \centering
     \includegraphics[width= .9 \textwidth]{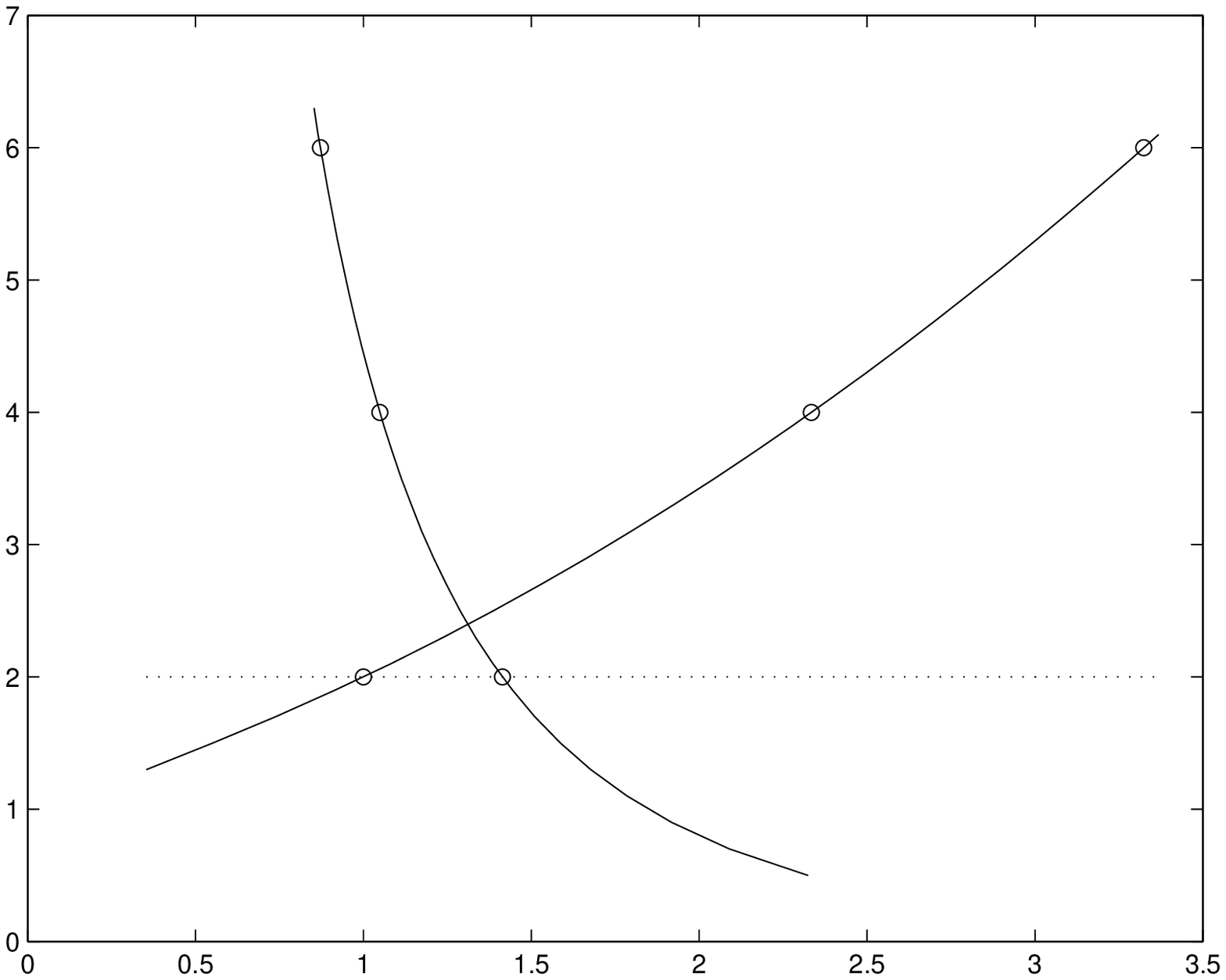}%
  \caption{Constructing $\, \alpha^+_1$.}
  \label{controcontro}
   \end{minipage}
\end{figure}
This construction could be achieved also by exploiting Lemmas \ref{sette} and \ref{otto} and defining $\, \alpha_1^+\,$ and $\,\alpha_1^-\,$ by solving the equations
$$\alpha^-(s) = \ove \alpha^+(s),\qquad \alpha^+(s) = \ove \alpha^-(s).$$
The solvability of the latter equations is evident from Lemmas \ref{sette} and \ref{otto} and the the corresponding common value of the functions is exactly $\, \alpha^-_1\,$ and $\, \alpha^+_1$, respectively. This construction is depicted in Figures \ref{contro}-\ref{controcontro} for the case $\, N=1$. The figures are outputs of actual computations and have been produced with {\small MATLAB}$^{\circledR}$. The dots correspond to the case of heat polynomials of even degree 2, 4, and 6 and have been graphically superimposed on the computed curve in order to validate the numerical results (in particular, the curves are not interpolated from the dots). Let us stress that the figures confirm once again that $\, 0 < \alpha_1^- < 2 < \alpha^+_1$.

\subsection{Construction of solutions changing sign $\, k \,$ times}

The argument developed above can be extended in order to construct solutions changing sign exactly $\, k \,$ times. 

Owing to the above discussion, for all $\, \alpha > \alpha^+_1 \,$ the solution $\, f\,$ to \eqref{S} with $\,f(0) >0 \,$ changes sign at least twice. In particular, one readily has that 
$$\lim_{\alpha \searrow \alpha^+_1}  s^{+,2}_\alpha =+\infty,$$
by construction. Hence, we readily find $\, \alpha > \alpha^+_1\,$ such that 
$$  s^{+,2}_\alpha > \ove s^{-,1}_\alpha.$$
Now, exploiting Lemma \ref{zz}, one may increase $\, \alpha \,$ up to the unique homogeneity $\, \alpha^{+}_2\,$ such that 
$$  s^{+,2}_{\alpha^+_2} = \ove s^{-,1}_{\alpha^+_2}$$
(note that $\, \alpha_1^+<\alpha^+_2$). Equivalently, one could define $\, s \mapsto \alpha^{+,2}(s)\,$ as the inverse of $\, \alpha \mapsto s^{+,2}_\alpha\,$ (which is well-defined due to Lemma \ref{zz}) and then solve.
\begin{equation}
  \label{spirit}
  \alpha^{+,2}(s)= \ove \alpha^+(s).
\end{equation}
The latter equation admits the unique solution $\, \alpha^+_2$, see Lemmas \ref{otto} and \ref{zz}. Figures \ref{contro2} below shows the (computed) outcome of the above described construction in the case $\, N=1$. In particular, in analogy with Lemma \ref{cosa}, a self-similar profile $\, f \,$ with homogeneity $\, \alpha > 1 \,$ and $\, f(0) >0\,$ has exactly two sign changes iff $\, \alpha^+_1 < \alpha \leq \alpha^+_2$. Moreover, among all these self-similar profiles, the only one that grows as $\, s^{\alpha}\,$ for large $\,s \,$ is the one corresponding to $\, \alpha= \alpha^+_2$.

The argument can be iterated for determining a strictly increasing sequence $\, \alpha^+_k\,$ such that all self-similar profiles $\, f \,$ with $\, f(0)>0\,$ changing sign exactly $\, k \,$ times ($ k >2 $) are given by the homogeneities $\, \alpha \,$ such that $\, \alpha^+_{k-1}< \alpha \leq \alpha^+_{k}$. Finally, $\, \alpha^+_k\,$ can be easily proved to be unbounded by contradiction. Hence, Theorem \ref{nuovomain}.i-iii follows.
Finally, it should be clear that the positivity assumption $\, f (0)>0 \,$ plays no specific role and the above discussion could be repeated for self-similar profiles which are negative in $\, 0 \,$ as well. This will eventually give rise to the corresponding sequence $\, \alpha^-_k\,$ and Figure \ref{controcontro2} reports the construction of $\, \alpha^-_2\,$ in the same spirit of \eqref{spirit}.

Let us close this subsection by remarking that Figures \ref{contro2} and \ref{controcontro2} show in particular that the ordering $\, \alpha^-_k < \alpha^+_k \,$ which holds at level $\, k =1 \,$ is however not to be expected in general (recall however that some ordering property holds for all $\, k \,$ holds by Theorem \ref{nuovomain}.v).

Finally, we shall explicitly observe that whenever $\, \alpha = \alpha^\pm_k \,$ for some $\, k$, we have that
\begin{equation}
  \label{ident}
  s^{\pm,j}_\alpha = \ove  s^{\pm,j}_\alpha \quad \text{for} \ \ j=1, \dots,k.
\end{equation}
\begin{figure}
   \begin{minipage}[t]{0.5\linewidth}
      \centering
      \includegraphics[width= .9 \textwidth]{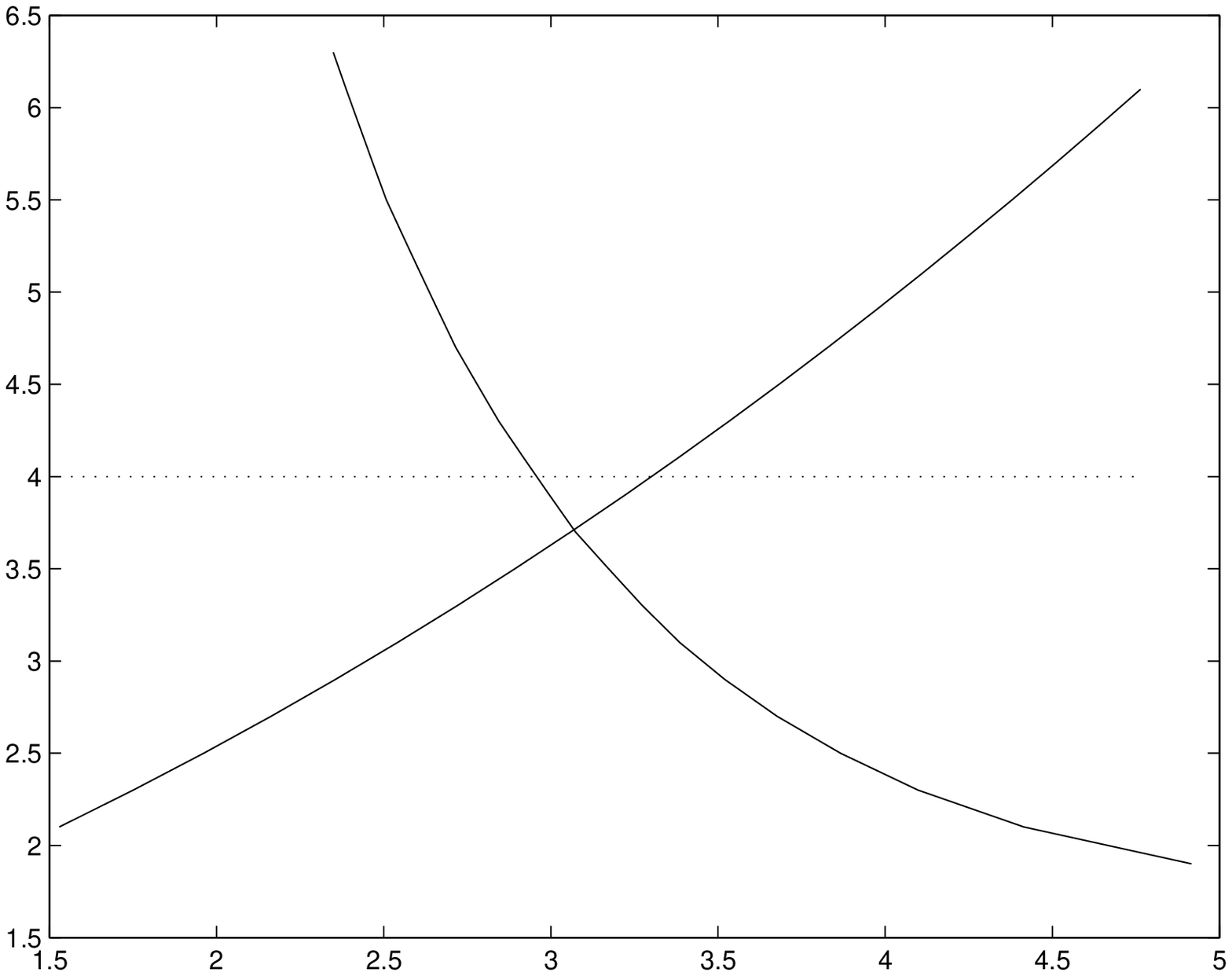}%
  \caption{Constructing $\, \alpha^+_2$.}
  \label{controcontro2}
   \end{minipage}%
   \begin{minipage}[t]{0.5\linewidth} 
     \centering
     \includegraphics[width= .9 \textwidth]{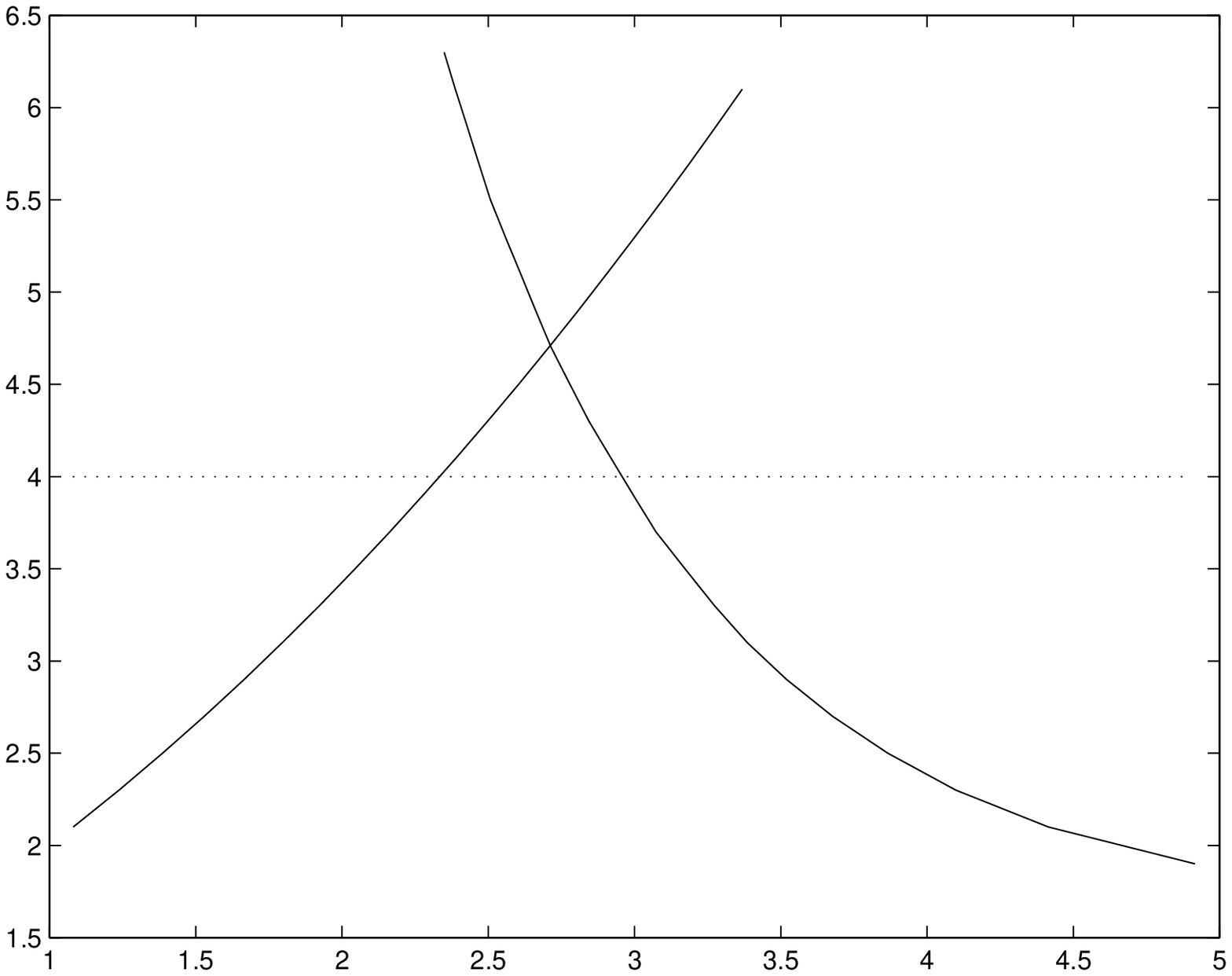}%
  \caption{Constructing $\, \alpha^-_2$.}
  \label{contro2}
   \end{minipage}
\end{figure}

\subsection{Proof of $\,\alpha^-_k< \alpha^+_{k+1}<\alpha^-_{k+2}$}
This fact is again to be obtained by exploiting Lemma \ref{zz} and explicitly matching shooting semi-profiles by triggering the corresponding homogeneity. In particular, by recalling \eqref{ident} we readily check that
$$  s^{+,1}_{\alpha^+_{k+1}} = \ove s^{+,1}_{\alpha^+_{k+1}}<  s^{-,2}_{\alpha^+_{k+1}} ,$$
where the inequality is a straightforward consequence of the minimality in Lemma \ref{polp}. In particular, by letting $\, f^\pm \,$ be the self-similar profiles corresponding to $\, \alpha^+_{k+1}\,$ with $\, f^\pm(0)=\pm 1$, one readily checks that
$$R^+\left(f^-,\left(s^{-,1}_{\alpha^+_{k+1}},s^{-,2}_{\alpha^+_{k+1}} \right)\right) = \alpha^+_{k+1} = R^+\left(f^+,\left(0,s^{+,1}_{\alpha^+_{k+1}} \right)\right),$$
which clearly leads to a contradiction if $\,  s^{+,1}_{\alpha^+_{k+1}} \geq  s^{-,2}_{\alpha^+_{k+1}}$, namely if
$$\left(s^{-,1}_{\alpha^+_{k+1}},s^{-,2}_{\alpha^+_{k+1}} \right)\subset\left(0,s^{+,1}_{\alpha^+_{k+1}} \right).$$

Hence, in order to find a self-similar profile $\, f \,$ which is negative in $\, 0 \,$ and changes sign $\, k+2 \,$ times, we shall increase the corresponding homogeneity above $\, \alpha^+_{k+1}\,$ and determine the unique value $\, \alpha^-_{k+2}\,$ such that
$$ \ove s^{+,1}_{\alpha^-_{k+2}}= s^{-,2}_{\alpha^-_{k+2}} .$$
Conversely, we readily check from \eqref{cambio} and \eqref{ident} that
$$  s^{-,1}_{\alpha^+_{k+1}} = \frac{s^{+,1}_{\alpha^+_{k+1}}}{\sqrt2}< s^{+,1}_{\alpha^+_{k+1}} <  s^{+,2}_{\alpha^+_{k+1}} = \ove s^{+,2}_{\alpha^+_{k+1}},$$
so that, by decreasing the homogeneity below $\, \alpha^+_{k+1}\,$ we eventually find a unique value $\, \alpha_k^-\,$ such that 
$$s^{-,1}_{ \alpha_k^-} = \ove s^{+,2}_{\alpha_k^-}.$$
The latter corresponds to the homogeneity of a self-similar profile which is negative in zero and changes sign $\, k \,$ times.

\section{An Appell transform argument}\label{Appell}
\setcounter{equation}{0}

Let $\, \phi \,$ be the (rescaled) fundamental solution of the heat equation and $\, w \,$ be a caloric function for $\, t <0$. Then, its Appell transform  \cite{Widder63}
$$v(x,t) := \phi(x,t)w\left(\disp\frac{x}{t}, -\disp\frac{1}{t}\right) \qquad\text{for} \ \  t>0$$
is again caloric. Moreover, the function $\, v \, $ is $\, -(N+\alpha)$-homogeneous iff $\, w \,$ is $\, \alpha$-homogeneous. Indeed, by letting
$$w(x,-t)=t^{\alpha/2}f(|x|/\sqrt{t}), \quad \phi(x,t)= t^{-N/2}\psi(|x|/\sqrt{t})\qquad \text{for} \ \ t>0,$$
with $\,\psi(r):= \exp(-r^2/4)\,$ for $\, r>0$, we readily check that
\begin{equation}\label{appell}
v(x,t)= t^{-(N+\alpha)/2}\psi(r)f(r)=:t^{-(N+\alpha)/2} g(r) \qquad \text{for} \ \ r=|x|/\sqrt{t}.
\end{equation}

We aim to exploit this construction in order to reformulate our results for self-similar solutions to \eqref{eq2} of negative homogeneity. The latter shall be defined for $\, t > 0\,$ and possibly develop a singularity as $\, t \to 0^+$. To this end, we suitably adapt Appell's transformation to the present situation. In particular, let us assume $\, w \,$ to be a radial self-similar solution to \eqref{eq2} for $\, t<0\,$ and denote by $\, \alpha >0 \,$ its homogeneity. Owing to Theorem \ref{nuovomain} and Lemma \ref{cosa}, we readily have that $\, w(\cdot, -1) \,$ changes sign exactly $\, k \,$ times. Hence, the corresponding self-similar profile $\, f_w \,$ has exactly $\, k \,$ zeros, namely $\, s_1, \dots, s_k$. Let $\, \psi_w : [0,+\infty) \rightarrow \Rz \,$ be defined as
\begin{equation}\label{psi} 
\psi_w(r)= \exp\left(-\disp\int_0^r \disp\frac{\lambda(f_w(s))s}{2}ds\right)\qquad \text{for} \ \ r \geq 0,
\end{equation}
where 
$$ \lambda (y) :=  1/\gamma(y) = 4/(3 + \text{sign}\,(y))\quad \forall y \in \Rz.$$
 Here, the function $\, \psi_w \,$ plays the role of the self-similar profile of the heat kernel $\, \phi_w(x,t):= t^{-N/2}\psi_w(|x|/\sqrt{t})\, $ whose discontinuous coefficients are driven by (the sign of) $\, w $. 

In analogy with \eqref{appell}, we shall check that the position
$$g:= \psi_w f_w,$$
gives actually rise to a solution to \eqref{eq2}. To this aim, let us explicitly observe that $\, g \in C^{1,1}[0,+\infty)\,$ since $\, f_w \,$ vanishes where $\, \phi_w \,$ is non-differentiable. Namely, $\, g \,$ is exactly as regular as $\, f_w$. 

Let us now observe that
\begin{gather}
\psi_w''(r) + \left(\disp\frac{N-1}{r} + \frac{\lambda(f_w(r))r}{2}\right)\psi_w'(r) + \disp\frac{\lambda(f_w(r))N}{2}\psi_w(r) =0 \nonumber\\
 \text{for} \ \ 0< s \not = s_1, \dots, s_k.\label{pp}
\end{gather}
More precisely, we stress that the latter equation is fulfilled in the distributional sense in all of $\, (0,+\infty)\, $ by including the measure
$$-\frac{\text{sign}(f_w(0))}{2} \sum_{i=1}^k (-1)^i\delta_{s_i} f'_w(s_i)s_i \psi_w(s_i)$$
(with obvious notation for the Dirac masses) in the right hand side. We shall now multiply equation \eqref{pp} by $\, f_w \,$ and add it to the corresponding equation in \eqref{S} multiplied by $\, \psi_w$. By exploiting the strict positivity of $\, \psi_w \,$ we directly have that $\, \text{sign}( g)= \text{sign}( f_w)$. Moreover, by making use of the fact that (see \eqref{psi})
$$ 2 \psi_w'(r) f_w'(r) = - \lambda (f_w(r)) r \psi_w'(r) f_w'(r) \quad \forall r >0,$$
and that $\, \delta_{s_i} f(\cdot)=0 \,$ for $\, i =1, \dots ,k$,
one computes that
\begin{eqnarray}
 \left\{
\begin{array}{l}
g(0)=f(0), \quad g'(0)=0,\quad\\
\disp  g''(r) + \left( \frac{N-1}{r}+r \right)g'(r) + (N+\alpha) g(r)=0 \quad \text{where} \ \ g <0, \ \ r>0 \\
\disp g''(r) + \left( \frac{N-1}{r}+\frac{r}{2} \right)g'(r) + \frac{N+ \alpha }{2} g(r)=0\quad \text{where} \ \ g>0, \ \ r>0 .
\end{array}
\right.\nonumber
\end{eqnarray}
Namely, the function 
$$v(x,t):=\phi_w(x,t)w\left(\disp\frac{x}{t},-\disp\frac{1}{t}\right) = t^{-(N+\alpha)/2}g\left(\frac{|x|}{\sqrt{t}}\right)$$
is a radial self-similar solution to \eqref{eq2} for $\, t>0$. Moreover, the bound \eqref{bound} is readily transformed into
\begin{equation}
 r \mapsto e^{\ell r^2/4}r^{-\alpha}g(r) \quad \text{is bounded},\label{bound2}
\end{equation}
where $\, \ell:= \lim_{r\to +\infty }\lambda(g(r))\,$ (which obviously exists). In particular, $\, v \,$ is such that
$$\lim_{t \searrow 0} v(x,t)= |x|^{-(N+\alpha)}\lim_{r \to +\infty} r^{N+\alpha} g(r) =0 \qquad \text{for} \ \ x\not =0.$$
More precisely, we have obtained that
\begin{equation}\label{bound3} |x|^{N+\alpha}\left(\disp\frac{|x|}{\sqrt{t}}\right)^{-(N+2\alpha)} e^{\ell |x|^2/4t} v(x,t) \quad \text{is bounded}.
\end{equation}

Since $\, \psi_w \,$ is strictly positive, it is a standard matter to check that, starting from a radial self-similar solution $\, v\,$ to \eqref{eq2} for $\, t >0 \,$ with homogeneity $\, -(N+\alpha)$, the position
\begin{equation}\label{inverse}
 w(x,t):= \phi^{-1}_v\left(\frac{x}{t}, -\frac{1}{t}\right)v\left(\frac{x}{t}, -\frac{1}{t}\right)\qquad\text{for} \  \ t<0
\end{equation}
(with obvious notation for $\, \phi_v$) gives rise to a radial self-similar solution to \eqref{eq2} for $\, t<0$. The latter relation can be rephrased for the corresponding self-similar profiles as $\, f = \psi_v^{-1} g$.
Moreover, the bound \eqref{bound3} obviously entails \eqref{bound}. It is finally straightforward to check that the above described modified Appell's transformation and its inverse \eqref{inverse} compose to the identity. 

Finally, by applying the above described transformation, we can recast our characterization result for the case of radial self-similar solutions of \eqref{eq2} for $\, t > 0\,$ with negative homogeneity.

\begin{theorem}\label{nuovomain2}
The sequences $\,\beta_k^+=-N-\alpha^+_k\,$ and $\, \beta_k^- =-N-\alpha^-_k\,$ are such that
\begin{enumerate}
\item[i)] for all $\, \beta <-N$, the unique $\, \beta$-homogeneous solution $\,v \,$ with $\, v(0,1)=\pm 1\,$ fulfills \eqref{bound3} iff $\, \beta = \beta_k^\pm \,$ for some $\,k \in \Nz$,
\item[ii)] for all $\, \beta <0$, the self-similar profile does not change sign iff
$\, -N\leq \beta<0$, changes sign once iff $\, \beta_1^\pm  \leq \beta < -N$, and changes sign exactly $\, k \,$ times ($k \geq 2$) iff $\, \beta_{k}^\pm  \leq \beta < \beta_{k-1}^\pm$.
\item[iii)] $\, \beta^\pm_k \,$ are strictly decreasing and unbounded,
\item[iv)] $\, \beta^+_1<  -N - 2 < \beta^-_1< -N$,
\item[v)] $\,\beta^-_{k+2}<\beta^+_{k+1}<\beta^-_k \ $ for all $\,  k \in \Nz$.
  \end{enumerate}
\end{theorem}

\begin{proof}[Sketch of the proof]
  We shall not report here a detailed proof. Indeed, the latter can be readily obtained by systematically exploiting the above described transformation and Theorem \ref{nuovomain} and Lemma \ref{cosa}. Let us just comment on Theorem \ref{nuovomain2}.ii. Namely, we aim to prove that, given $\, -N \leq \beta <0$, and $\, g(0)=1$, for instance (the case $\, g(0)=-1\,$ being completely analogous) the corresponding solution to 
\begin{equation}\label{usa0}
\disp g''(r) + \left( \frac{N-1}{r}+\frac{r}{2} \right)g'(r) - \frac{\beta }{2} g(r)=0,\qquad g'(0)=0 
\end{equation}
does not change sign. By contradiction let $\, g(r)=0\,$ for some $\, r >0$. Hence, the transform argument and the eigenvalue characterization given by Lemma \ref{polp} entail that
$\, \beta = - N - \alpha^+(r) < -N$, contradicting the fact that $\, -N < \beta$.

More precisely, one can prove that $\, g(r)> \exp(-r^2/4)\,$ for all $\, r >0$. Indeed, we readily compute that $\, p(r):= g(r) - \exp(-r^2/4)\,$ fulfills
\begin{equation}\label{usa}
\disp p''(r) + \left( \frac{N-1}{r}+\frac{r}{2} \right)p'(r) + \frac{N }{2} p(r)= \frac{N+\beta}{2} g(r) >0\qquad \text{where} \ \ g >0.
\end{equation}
Moreover, we easily check that $\, p' >0 \,$ in a right neighborhood of $\, 0\,$ (arguing for instance as in the proof of Lemma \ref{zero}).
Assume by contradiction that there exists $\, r >0 \,$ such that $\, p(r)=0\,$ and $\, p >0 \,$ on $\, (0,r)$. Then, relation \eqref{usa} entails that
\begin{gather}
\left(\disp \int_0^r s^{N-1} e^{s^2/4} (p'(s))^2ds\right)\left(\disp \int_0^r s^{N-1} e^{s^2/4} p(s)^2ds\right)^{-1} <- \frac{N}{2}.\label{pl}
 \end{gather}
On the other hand, since $\, r \mapsto \exp(-r^2/4)\,$ never vanishes and is the only solution to \eqref{usa0} for $\, \beta = -N$, by arguing along the same lines of Lemma \ref{polp} one readily has that 
\begin{gather}
 \left(\disp \int_0^r s^{N-1} e^{s^2/4} (q'(s))^2ds\right)\left(\disp \int_0^r s^{N-1} e^{s^2/4} q(s)^2ds\right)^{-1} \geq -\frac{N}{2} \nonumber\\ \forall q \in L^2((0,r),\mu) \  \ \text{such that} \ \ q \not = 0, \ q' \in  \ L^2((0,r),\mu), \ q(r)=0,\nonumber
\end{gather}
where $\, \mu \,$ is the absolutely continuous measure $\,d\mu(s) = s^{N-1} e^{s^2/4} ds$.
This contradicts \eqref{pl} and the assertion follows.
\end{proof}

\subsection{The result by Kamin, Peletier, \& Vazquez} As already mentioned, the analysis of the radial self-similar solvability of \eqref{eq} has been addressed in \cite{Kamin-Peletier-Vazquez91} for $\, t>0 $. In particular, {\sc Kamin, Peletier, \& Vazquez} focus on self-similar solutions of the form
$$u(x,t)= t^{\eta/2} m(r) \qquad \text{for} \ \ r=|x|/\sqrt{t}, \ \eta < 0 $$
and develop a full analysis on their relevance in determining  the long-term behavior of solutions to \eqref{eq} with exponentially decaying, continuous, and non-negative initial data.

Among these results, they prove in particular that there exists a unique $\, \eta^-_1 \,$ such that the corresponding self-similar solution $\, u \,$ is positive and $\, r^{-\eta^-_1}m(r) \to 0 \,$ as $\, r \to 0\,$ \cite[Thms. 2.1-2.2]{Kamin-Peletier-Vazquez91}. Moreover, they directly check that 
$$ \eta^-_1 < - N +2 $$ 
and it is such that $\, u_t(\cdot,1) \,$ changes sign exactly once being negative inside some given ball and positive outside \cite[Lemma 2.12]{Kamin-Peletier-Vazquez91}.

Our analysis is independent from the former since, clearly, the time derivative of a self-similar solution is not self-similar in general. On the other hand, the above-mentioned results are strongly reminiscent of our situation where we obtain that, letting $\, w \,$ be a negatively-homogeneous radial solution to \eqref{eq2} for $\, t>0$, the function $\, w(\cdot,1) \,$ changes sign exactly once in correspondence to some specific homogeneity range. Moreover, within the latter range we can select a specific value 
$$ \beta^-_1 < -N $$ 
such that $\, w(x, \cdot) \,$ behaves nicely as $\, t \searrow 0$. Hence, at least at this formal level, our value $\, \beta^-_1 \,$ plays essentially the role of $\, \eta^-_1 -2 \,$ in the analysis of \cite{Kamin-Peletier-Vazquez91}.

\bibliographystyle{plain}
\def\cprime{$'$} \def\cprime{$'$}

\end{document}